\documentclass[12pt,A4paper]{article}

\newcommand{\X}{\mathcal{X}}
\newcommand{\Z}{\mathcal{Z}}

\newcommand{\Tcal}{\mathcal{T}}
\newcommand{\Ncal}{\mathcal{N}}
\newcommand{\Ecal}{\mathcal{E}}
\newcommand{\Ccal}{\mathcal{C}}
\usepackage{amsmath}

\usepackage{times}
\usepackage{bm}
\usepackage{natbib}
\usepackage{gensymb}
\usepackage[plain,noend]{algorithm2e}

\usepackage{amsmath}
\usepackage{amssymb}
\usepackage{mathrsfs}
\usepackage{amstext}
\usepackage{graphicx,epsf}
\usepackage{enumerate}
\usepackage{natbib}
\usepackage{color}
\usepackage{multirow}

\makeatletter
\renewcommand{\algocf@captiontext}[2]{#1\algocf@typo. \AlCapFnt{}#2} 
\def\@algocf@capt@plain{top}
\renewcommand{\algocf@makecaption}[2]{%
  \addtolength{\hsize}{\algomargin}%
  \sbox\@tempboxa{\algocf@captiontext{#1}{#2}}%
  \ifdim\wd\@tempboxa >\hsize
    \hskip .5\algomargin%
    \parbox[t]{\hsize}{\algocf@captiontext{#1}{#2}}
  \else%
    \global\@minipagefalse%
    \hbox to\hsize{\box\@tempboxa}
  \fi%
  \addtolength{\hsize}{-\algomargin}%
}
\makeatother

\begin{document}



\markboth{K. GARSIDE ET AL}{Event History Analysis for Topological Data}

\begin{center}
\begin{Large}
{\bf  Event History and Topological Data Analysis}
\end{Large}

\vspace*{1cm}

KATHRYN GARSIDE, AIDA GJOKA, ROBIN HENDERSON, HOLLIE JOHNSON, IRINA MAKARENKO

School of Mathematics, Statistics and Physics, Newcastle University, NE1 7RU, U.K. 

\end{center}

\begin{abstract}
Persistent homology is used to track the appearance and disappearance of features as we move through a nested sequence of topological spaces. Equating the nested sequence to a filtration and the appearance and disappearance of features to events, we show that simple event history methods can be used for the analysis of topological data.  We propose  a version of the well known Nelson-Aalen cumulative hazard estimator for the comparison of topological features of random fields and for testing parametric assumptions.  We suggest a Cox proportional hazards approach for the analysis of embedded metric trees.  The Nelson-Aalen method is illustrated on globally distributed climate data and on neutral hydrogen distribution in the Milky Way.  The Cox method is use to compare vascular patterns in fundus images of the eyes of healthy and diabetic retinopathy patients.

\end{abstract}


\section{INTRODUCTION}

Although interest in topological data analysis has grown sharply in recent years, most of the  important developments have appeared outwith the statistical literature.  We suggest that a consequence is that the potential benefits of adapting some established statistical methods to topological data may not have been fully exploited.  In this paper we argue that event history methods are highly suitable for the analysis of data on certain topological features, that the use of such methods can increase the power and effectiveness of topological data analyses, and in turn that exploitation of topological features can bring added value to statistical analyses.   

\begin{figure}
\centerline{\includegraphics[height=5in,width=5in]{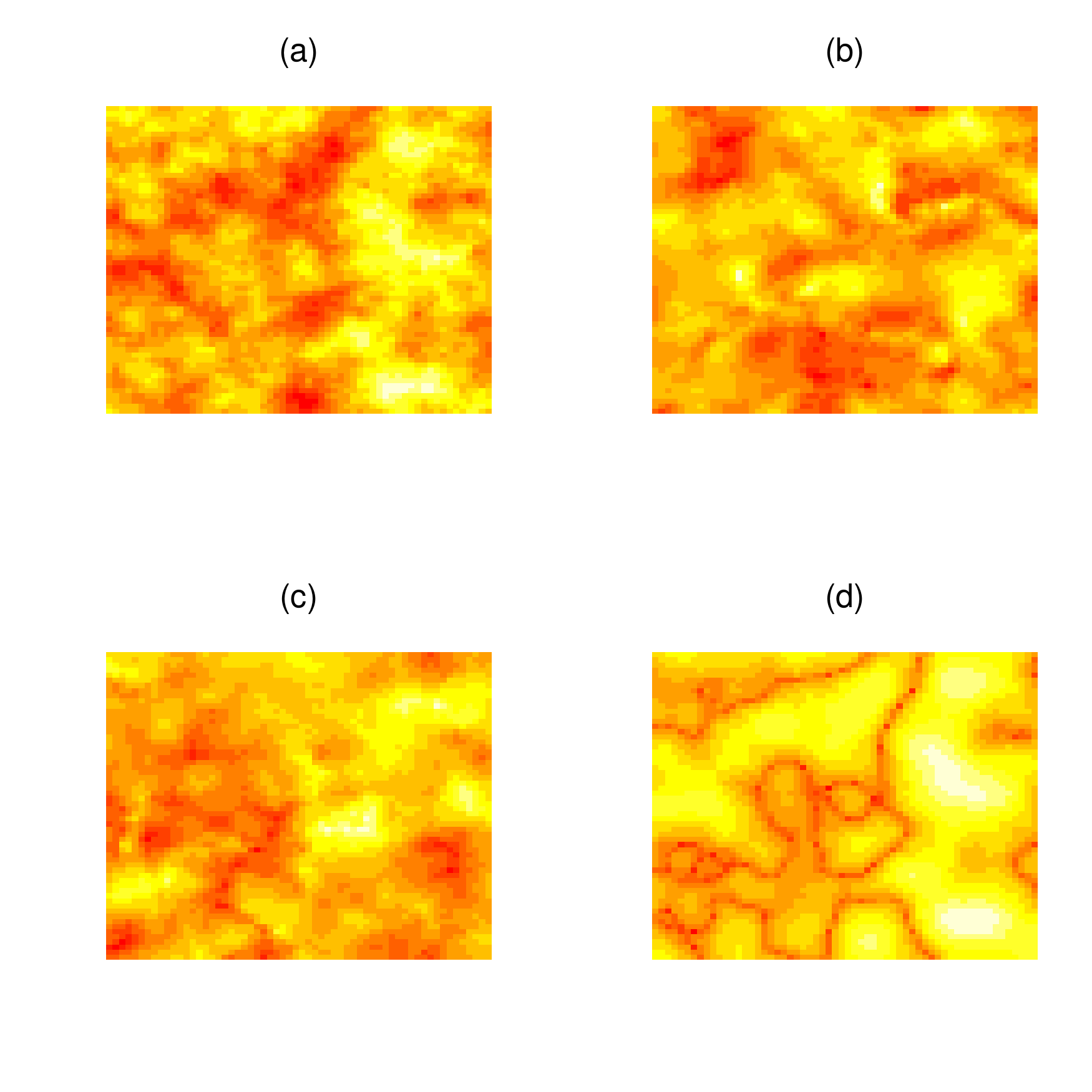}}
\caption{Four simulated random fields, each on a 60$\times$60 lattice.}
\label{fig:intro1}
\end{figure}

\begin{figure}
\centerline{\includegraphics[height=3in,width=3in]{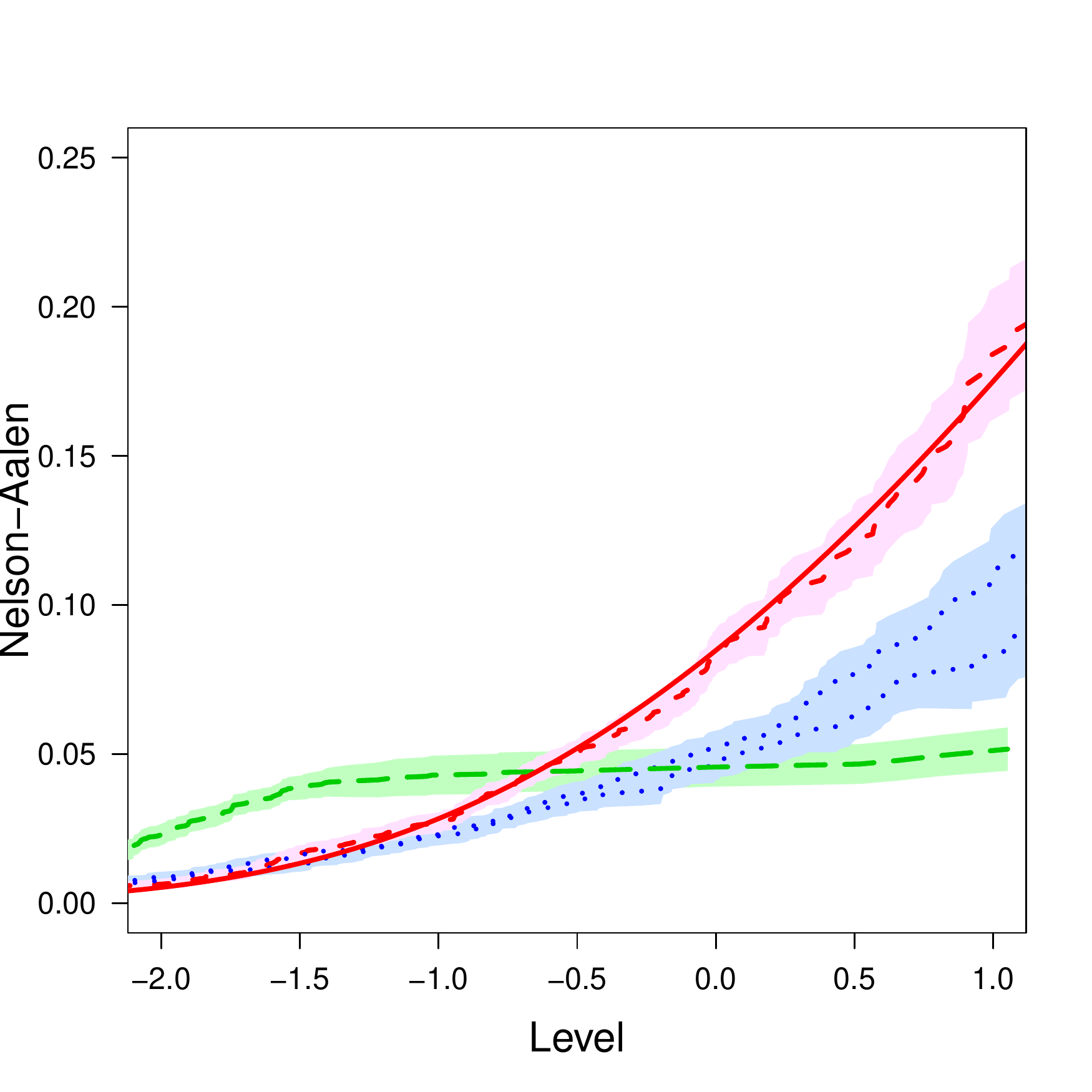}}
\caption{Nelson-Aalen plots for connected components for the random fields of Fig. \ref{fig:intro1}.  Short dashed lines correspond to panel (a), dotted  lines to panels (b) and (c), and long dashed lines to panel (d).  The shaded regions indicate $\pm$ one standard deviation.  The solid line shows the expected value for a Gaussian random field with the same correlation structure as the data in Fig. \ref{fig:intro1}.}
\label{fig:intro2}
\end{figure}

We use Fig. \ref{fig:intro1} and  Fig. \ref{fig:intro2} to illustrate, at the cost of preempting Section 3, where  they will be explained more fully. Figure \ref{fig:intro1} shows four simulated stationary and isotropic random fields, each on a $60\times 60$ discrete lattice. They are simulated to have correlation functions that are indistinguishable for practical purposes. All four have standard N(0,1) marginal distributions. 
 Panel (a) is a simulation from a Gaussian random field, panels (b) and (c) are independent realisations of one particular non-Gaussian random field, and panel (d) is a realisation of a second non-Gaussian random field.  Panel (d) is  visually quite different from the others, but we see no obvious differences between (b) or (c) and the Gaussian random field (a).  Figure \ref{fig:intro2} however shows simple  Nelson-Aalen plots based on topological features of these data, to be described in Section 3. Now the differences are stark.    Given that all four have the same marginal and correlation properties, we know of no other method that can draw out higher order differences so readily.  Further, panel (a) is  seen to be consistent with a Gaussian random field while the other panels are not.

Topological data analysis is in part founded on persistent homology. Important early contributions include   \cite{frosini1992measuring},   \cite{robins1999towards} and, particularly, \cite{Edelsbrunner2002}, who introduced the idea of persistence of features along filtrations.  However, the earliest appearance of the expression ``topological data analysis'' seems not to have been until work by \cite{de2004topological} and \cite{bremer2004topological}. Later, \cite{carlsson2009topology} was key in the popularisation of topological data analysis, demonstrating how ideas from topology could solve many of the issues faced when applying geometric methods to complex data.  Other  developments are described by \cite{perea2019}.


In Section 2 we give an introduction to  the key  constructs and concepts underpinning persistent homology and topological data analysis, aimed at a statistical reader.  In Section 3 we develop the idea of Nelson-Aalen plots for  topological features of random fields, as already seen in our introductory example.  We apply the methods in two application areas, one being properties of climate model residuals on a sphere, and the other the distribution of neutral atomic hydrogen \textup{H}\,{\sc i} in the interstellar medium. In Section 4 we consider topological methods for embedded metric trees and propose Cox regression methods as a simple way to deal with potentially censored tree data.  This method is applied to vascular structures in fundus images of the human eye, where it has potential to be an effective screening tool for the early detection of diseases.

\section{TOPOLOGICAL DATA ANALYSIS}

Topology is the study of features in a  space that are invariant under continuous deformation.  We will concentrate on topological data analysis methods that are based on persistent homology, which is the most popular variant of topological data analysis.  Detailed descriptions  are provided by, inter alia,  \cite{Edelsbrunner2002}, \cite{edelsbrunner2008persistent}, \cite{carlsson2009topology} and \cite{Edelsbrunner2012PersistentHT}, with an excellent computational description given by \cite{Otter2017}.  Other variants of topological data analysis are summarised by \cite{wass18}.

We assume we have a metric space and inherited distance function. Dealing with arbitrary features in a space is complicated, and so a starting point is to give a representation in a common form.  One such is a simplicial complex, which is a collection of vertices, edges, triangles, tetrahedra and higher dimensional equivalents.  Thus the letters T, D, and A when written on a page could be represented by four vertices connected by three edges (T), a triangle (D) and a triangle with two extra external edges to two extra vertices (A). 
Simplicial complexes are not unique but they facilitate the use of linear algebra for the computation of topological invariants, which are features that do not change when the underyling space is bent, stretched or otherwise deformed, but without  gluing or tearing.  

The most important topological invariants are the homology groups, which are formally algebraic constructs, one per natural number, that collect together invariant features at each dimension. The ranks of these groups are called Betti numbers. The zero-dimension homology group $H_0(\X)$ of a simplicial complex $\X$ consists of the disjoint elements, often referred to as connected components.   For example, if the letters T, D, and A are written on a page then there are three connected components and the zero-order Betti number is $\beta_0=3$. In three dimensions, if we take for example letter balloons of T, D and A,  there are still three connected components.

The next homology group $H_1(\X)$ consists of holes constructed from one-dimensional edges. These are loops or cycles that cannot be shrunk to a single point. For T, D, and A on a page there are two, so $\beta_1=2$.  For  T, D, and A as letter  balloons  we have $\beta_1=4$.  This is because the D and the A are each homologous to a torus.  The hole through the centre of the torus is obvious: the other hole arises because a loop around the circle of revolution that forms the torus cannot be shrunk to a point without leaving the surface.

The second-dimension homology group $H_2(\X)$ consists of two-dimensional cycles which enclose an empty three-dimensional space, often called a void. For T, D, and A in balloon form there are three, with one per connected component.  For T, D, and A on a page there are none.

\begin{figure}
\centerline{\includegraphics[height=5in,width=5in]{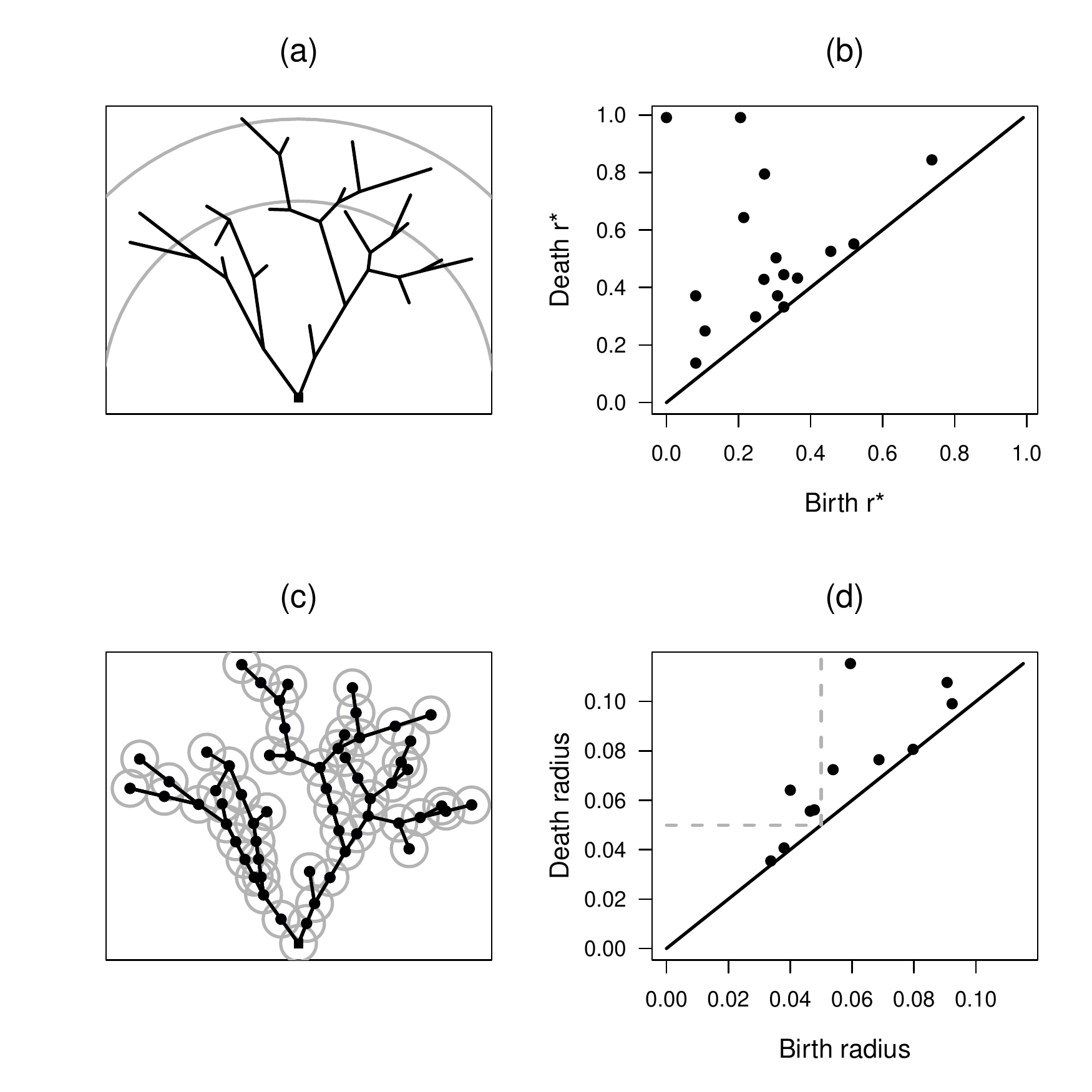}}
\caption{Illustrations of persistent homology: (a) embedded metric tree filtered by radius from root; (b) corresponding persistence diagram for connected components; (c) point cloud filtered by ball diameter; (d) corresponding persistence diagram for holes.  
}
\label{fig:trees}
\end{figure}

We now turn to persistent homology, which begins with a filtration of a space or simplicial complex $\X$. This is a nested sequence of subsets  
\begin{equation}
	\emptyset = \X_0 \subseteq \X_1 \subseteq \cdots \subseteq \X_d \cdots \subseteq \X_n  = \X,
\label{eqn:filtration}
\end{equation}
indexed by a parameter $d$ that can be discrete or continuous.  Each subset 
is itself a valid simplicial complex, with associated homology groups. Panel (a) of Fig. \ref{fig:trees} illustrates.  It shows a rooted metric tree embedded in two dimensions, and is a simplified version of the metric trees we will consider  in Section 4. There are no cycles and so no holes. In this example the feature of interest  is the connectivity of the tree as we move towards the root, filtering inward by radial distance. So $\X_d$ consists of the subset of $\X$ that is outside a circle of radius $r=1/d$ from the root.  The plot shows two examples.  The complex outside the outer circle consists of just one connected component.  The complex outside the inner circle consists of nine connected components.

In persistent homology the appearance and disappearance of features is tracked as we move through the filtration.  A feature is said to be born when it appears and to die when it disappears or merges with another feature.  Thus if we start with $r=\infty$ in Fig. \ref{fig:trees} there are no features.  This continues as $r$ is decreased  until it reaches the point on the tree that is furthest from the root, at radius $r_{\rm max}$ say.  The first component is born at this level.  As $r$ is decreased further, additional components are born at vertices that have only inward facing edges.  Existing components merge at vertices with two outward facing edges. In  these cases, by convention the component that was born first is assumed to continue, and the component that was born second is said to die at this radius. The process continues until the root is reached and the tree consists of a single connected component.

The birth and death levels of components are often represented by a persistence diagram, as in panel (b)  of Fig. \ref{fig:trees}.  This is a plot of birth times and death times of features as we move through the filtration.  In the example we have plotted against  $r_{\rm max}-r$ so that the first component to be born is on the left of the horizontal axis. Points near the diagonal are short-lived, points further from the diagonal are longer-lived.  In the example we have two subtrees with common root, which merge only at $r=0$ and so there are two points with the same maximum death time.

A common way to represent a topological space is through a point cloud, as in 
 panel (c)  of Fig. \ref{fig:trees}.  In these cases a filtration is usually defined by growing balls of increasing diameter, centred on the points. There are different ways of forming simplicial complexes in these cases, the most common of which is the Vietoris-Rips complex, which defines a subset of points to be a  simplex in $\X_d$ when all pairwise distances between points are less than $2d$.  The persistence diagram in panel (d) relates to the birth and death of holes for this example.  The grey broken lines indicate the ball radius shown in (c).  We see that a medium sized hole is about to be born towards the upper right and this will have modest persistence.  A large hole in the centre will be born when the ball radius is increased a little more.  This will be the last hole to die and is represented by the highest point in panel (d).

\begin{figure}
\centerline{\includegraphics[height=5in,width=5in]{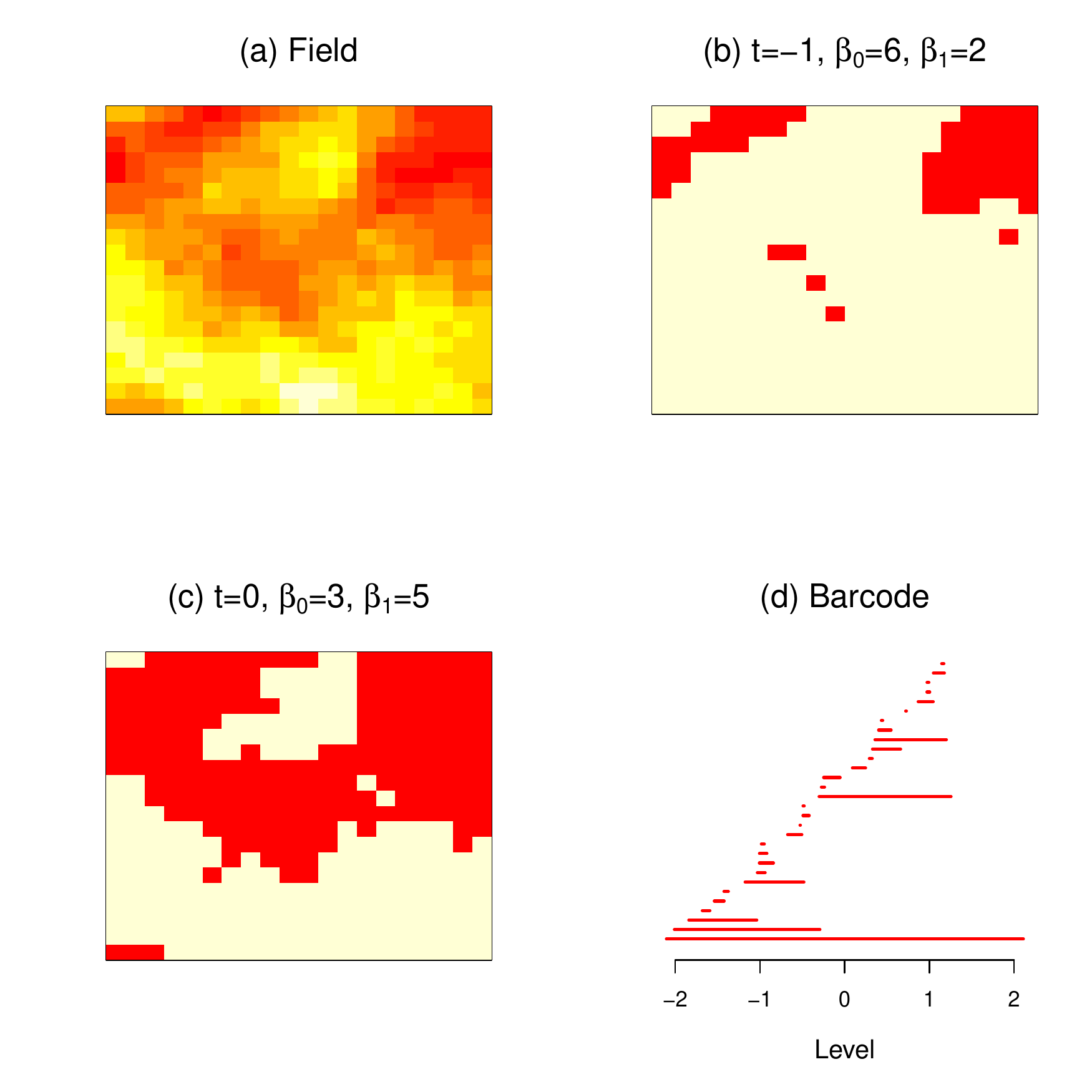}}
\caption{A simulated Gaussian random field  together with level sets  $\X_{-1}$ and $\X_{0}$, and barcode for connected components.}
\label{fig:levsets}
\end{figure}

Topological data analysis is often used when there is a real-valued function $z_x$ associated with each point $x$ in a space $\X$.  A lower level set $\X_t$ is defined as
\begin{equation}
 \X_t= \{x: z_x\leq t, x \in \X \},
\label{eqn:levelset}
\end{equation}
and upper level sets can defined similarly.  A filtration is defined by increasing or decreasing the threshold $t$ and persistence of topological features in the level sets is tracked.  Figure \ref{fig:levsets} illustrates.  In panel (a) we show a two-dimensional image on a discrete $20 \times 20$   grid.  We assume that pixels that share an edge are connected but those that share only a vertex are not. Panels (b) and (c) show two level sets: pixels with image values less than  $t=-1$ and  less than $t=0$ respectively, marked as dark squares.  In panel (b) there are $\beta_0=6$ connected components.  There are also two regions in that panel that are not in the level set and are completely separated by connected components.  Hence $\beta_1=2$. The level set has expanded by $t=0$ in  (c) and the six previous components have merged, along with others, to form one large component.  One new component has appeared to the lower left, and a small single-pixel component is close by, so $\beta_0=3$.  Five holes are evident. 

The final panel of Fig. \ref{fig:levsets} shows a so-called barcode for connected components for the image in panel (a).  This is an alternative to a persistence diagram for the presentation of lifespans of topological features as we filter through the level sets.  There is one line per feature, starting at the birth level and ending at the death level.  In the example the lowest field value is in the dark region to the top right.  This corresponds to the first pixel to appear in the level set, which by convention is the last to die, giving the long line at the bottom of the barcode  in (d).  The second component to be born is at a local minimum toward the top-left of the image.  This component also has the second-longest lifetime until it merges with the first at around $t=-0.25$.   The third-longest lifetime corresponds to the cluster at the lower left of the figure, which is relatively late to be born but quite isolated from other components and hence persists for some time.

A variety of approaches for the analysis of persistence diagrams and barcodes have been proposed.  Wasserstein or bottleneck distances can be used to measure the difference between two persistence diagrams, but these measures are not calibrated to any interpretable scale. There remains no useful equivalents of, say, means or variances for persistence diagrams, which would be needed for inference.  \cite{mileyko11}, \cite{turner14} and  \cite{munch15} showed how to obtain a {F}r\'{e}chet  mean for a group of persistence diagrams, but unfortunately this is not unique. An alternative is to base inference on summaries such as the landscape functions of 
\cite{bubenik2015}, which essentially consider  nested profiles of the point patterns,  or on counts of points over subregions, possibly with smoothing to produce the persistent images of, for example,  \cite{adams2015}. Other proposals include a form of parametric bootstrap suggested by 
\cite{adler17} and the accumulated persistence suggestion of \cite{bisc19}. 
\cite{wass18} provides a recent review.  Our approach is quite different.

\section{NELSON-AALEN FOR RANDOM FIELDS}

\subsection{Method}

Let $\Z$ be a zero-mean real-valued random field defined on a finite discrete space $\X$.  We assume a discrete space because in the applications we have in mind the field will often be an image, and invariably pixellated.   Let 
 $z_x$ be the field value at location $x \in \X$. We assume $z_x$ is drawn from a continuous distribution with no jump discontinuities and that there are no ties.  We will require cor$(z_x,z_{x^\prime})<1$ for all $x \neq x^\prime$ in $\X$.

Our filtration corresponding to \eqref{eqn:filtration} will be the nested sequence
of sub-level sets (\ref{eqn:levelset}).   But instead of tracking the persistence of features, which underpins most published applications of topological data analysis, we will take a counting process approach based on the birth levels of new features in the filtration. As well as being computationally much more simple, the approach leads naturally to the use of methods that are familiar in event history analysis, with level playing the role of time  and with negative values allowed.  As features of interest we will concentrate here on connected components within $\X_t$, namely the  $0$-dimensional features. 

A component is born at location $x$ at level $t$ if  $z_x=t$ and $x$ corresponds to a local minimum of $\Z$.  To define a local minimum in the discrete space $\X$ we need to define the neighbours $(x)$ of each location $x$.  We are free to do this as we wish, though shared edges or nearest neighbours seem sensible for pixellated or lattice-based fields.

Let $N_x(t)$ count the number of components that are born at location $x$ up to and including level $t$, and let $Y_x(t)$ be the associated predictable at-risk indicator. This is an indicator of whether, given the evolution of level sets to just before level $t$, it is possible for location $x$ to hold a local minimum at $t$.  So
\begin{equation*}
N_x(t) = \left\{ \begin{array}{ll}
1\;\;\;\; & z_x\leq t,\;\; z_{(x)}> z_x 1_x,\\
\\
0& {\rm otherwise}, \end{array} \right. 
\;\;\;\;\;\;\
Y_x(t) = \left\{ \begin{array}{ll}
1\;\;\;\; & z_x \geq t,\;\; z_{(x)}> t 1_x,\\
\\
0 & {\rm otherwise}, \end{array} \right. 
\label{eqn:counts}
\end{equation*}
where $1_x$ is a unit vector of the same length as $z_{(x)}$. Define $N(t)=\sum_xN_x(t)$ and $Y(t)=\sum_xY_x(t)$. Now let
\begin{equation}
\hat{A}(t) = \int_{-\infty}^t \frac{dN(u)}{Y(u)}.
\label{eqn:nelaal}
\end{equation}
This is an analogue of the Nelson-Aalen estimator of the cumulative hazard for a survival time variable \citep{aal08}. If preferred, we can redefine the level scale $t$ to $t^\ast=\exp(t)$ so as to be a more familiar non-negative argument. Like the Nelson-Aalen estimator,  $\hat{A}(t)$ is defined 
for all $t$ and has a finite number of discontinuities at the observed birth times. It is non-negative and uniformly bounded for fixed $\X$.    It is a cadlag process, meaning continuous to the right with left limits. However,   the 
standard Nelson-Aalen variance estimator
\begin{equation}
 {\hat{\rm{var}}\{\hat{A}(t)}\} = \int_{-\infty}^t \frac{dN(u)}{Y^2(u)} 
\label{eqn:var}
\end{equation}
does not in general apply, for a variety of reasons.  One is possible non-independence between counting processes at different locations given non-independence of the underlying field values.  Another, related, reason is informative censoring:  knowledge that location $x$ falls out of the risk set at level $u$ provides information about other locations.  And a third reason is that we have not constructed a martingale with respect to the observed data filtration.

We propose such Nelson-Aalen estimators for topological features as tools for comparing  properties of fields beyond marginal and correlation structures, for comparing observed field data or model residuals with an assumed parametric model, and as a diagnostic for assumed or fitted correlation structures.  We allow either a single field $\Z$ or   independent and identically distributed replicates $\Z_1, \Z_2, \ldots,\Z_N$ defined on a common space $\X$.

\subsection{Asymptotics}

If $\X$ is fixed and there are $N$ independent replicates $\hat{A}_1(t), \ldots, \hat{A}_N(t)$, then as each $\hat{A}_i(t)$ is non-negative and bounded above, for each level $t$ the variance of $\hat{A}_i(t)$ is finite  and the classic central limit theorem applies.  Hence pointwise inference for large $N$ is straightforward, based on the sample mean $\bar{A}(t)$ and associated sample standard deviation. 

For simultaneous confidence bounds we need to be a little careful because of the discontinuities in $\hat{A}_i(t)$.  \cite{hahn77} for instance points out that there can be uniformly bounded cadlag processes for which no central limit theorem applies, depending upon properties of their jumps.  \cite{bloz94} provide sufficient conditions under which $\surd N\bar{A}(t)$ converges weakly to a Gaussian process with finite variance and smooth sample paths in $(0,\tau)$, for any arbitrarily large $\tau$ for which pr$\{Y(\tau)>0\}>0$. These conditions are verified in Supplementary Material.  Hence a functional central limit theorem applies.

Further, for fixed $\X$  we can show that $\bar{A}(t)$ is consistent in $N$ for 
\begin{equation}
  A(t) = \displaystyle\int_{-\infty}^t J(u)\frac{ \sum_x {\rm pr}(z_{(x)}> u 1_x  \mid z_x=u) f_x(u)}{\sum_x {\rm pr}(z_x>u, \;z_{(x)}> u 1_x)} du,
\label{eqn:limit}
\end{equation}
where $f_x(.)$ is the marginal density of $z_x$ and $J(u)=1$ if $E\{Y(u)\}>0$, zero otherwise. If $N=1$ and there is just a single realisation, the estimator $\hat{A}(t)$ is also consistent for (\ref{eqn:limit}), but this time as the  cardinality of $\X$ increases, provided  cor$(z_x, z_{x^\prime})$ approaches zero as the distance between $x$ and $x^\prime$ increases.

\subsection{Inference}

For inference, we need to distinguish between the single and multiple replication situations.  We will begin by assuming multiple replications: $\hat{A}_1(t), \ldots, \hat{A}_N(t)$.

We can consider the replicates $\hat{A}_i(t)$ as functional data.  It is helpful to discretise the time/level scale of the estimator (\ref{eqn:nelaal})  to a  grid of $M$ distinct levels $\tau_1=t_1,\ldots,t_M=\tau_2$.  This is a common assumption in functional data analysis and deals effectively with the discontinuities in our estimators, which with probability one will not occur at measurement points.  
There are then a huge array of analysis techniques immediately available: see for example the excellent review by  \cite{wang16} for general approaches, and 
\cite{chiou09} for a specific example of treating discretised hazards as functional data.  With our applications in mind, we will mention only the construction of simultaneous confidence bands.  This problem has had attention for functional data that are either sparse  \citep{yao05, ma12} or dense \citep{degras11, cao12}.  Given that our discretisation is by design, a dense functional data approach is appropriate.

Let $\hat{A}_{ij}$ be the Nelson-Aalen estimator for field $\Z_i$ at level $t_j$.  Let $\bar{A}_j$ and $\hat{\sigma}_j^2$ be the sample mean and sample variance of the estimators at $t_j$, and let $\hat{\rho}_{jk}$  be the estimated correlation between $\hat{A}_{.j}$ and $\hat{A}_{.k}$.   \cite{degras11} 
establishes that for fixed $\tau_1$ and $\tau_2$, as both $N$ and $M$ increase 
\begin{equation*}
\bar{A}_j\pm c_\gamma\hat{\sigma}_j/\surd N
\label{eqn:ci}
\end{equation*}
has asymptotic coverage $1-\alpha$, where $c_\alpha$ is the upper $\alpha$-quantile 
of the maximum absolute value over $(\tau_2-\tau_1)/\tau_2$ of a Gaussian process with standard margins and with correlation function equal to an appropriately scaled limit of $\hat{\rho}_{jk}$.  The threshold $c_\gamma$ needs to be obtained numerically:  software is available  in the R package {\tt SCBmeanfd} \citep{scb16}.

Turning to the case of a single replication only, we assume a parametric model $F_\theta$ is available for the random field $\Z$ and we  propose a parametric bootstrap approach.  Parametric and more general bootstrap methods for functional data have been studied empirically and theoretically by, for example, \cite{cuevas06}. We require 
a consistent estimator $\hat{\theta}$ of the parameter vector $\theta$. We then generate $B$ simulated realisations of $\Z$ on $\X$ from $F_{\hat{\theta}}$, and  obtain the bootstrap Nelson-Aalen estimators $\hat{A}_b(t)$ for $b=1,\ldots,B$.  Pointwise intervals that have  correct coverage asymptotically in $B$ and $\mid \X \mid$ can be obtained from the bootstrap quantiles or a Normal approximation justified by the central limit theorem. 

For simultaneous confidence bands we  suggest a Monte Carlo method used by \cite{crain11} that avoids explicit estimation of the correlation function.  We again discretise the original and bootstrap estimators $\hat{A}_b(t)$ to a dense grid of $M$ points over a fixed range of levels $(\tau_1, \tau_2)$. Let $\hat{A}_j$ be the Nelson-Aalen  estimate (\ref{eqn:nelaal}) at level $t_j$ obtained from the original data, and let $\hat{A}_{bj}$ be corresponding value from bootstrap replicate $b$. Let $\bar{A}_j$ and $\hat{\sigma}_j$ be the bootstrap mean and estimated standard deviation at  $t_j$. Define 
\[ G_b= \max_j\{ \mid \hat{A}_{bj}- \bar{A}_j \mid/\hat{\sigma}_j \},\]
and take $\hat{d}_\alpha$ to be the upper $\alpha$-quantile of the empirical distribution of $G_1, \ldots, G_B$.  Then for large $M$ an approximate simultaneous confidence band for $A(t)$ over  $(\tau_1, \tau_2)$ is
\[ \hat{A}_j\pm  \hat{d}_\alpha \hat{\sigma}_j \hspace*{1cm} j=1,\ldots, M.\]

\subsection{Simulation studies}

We consider random fields on $d\times d$ lattices and investigate three models. The first,
M1,  is a stationary and isotropic Gaussian random field with standard N(0,1) marginals and Mat\'{e}rn correlation function.  If locations $x$ and $x'$ are separated by Euclidean distance $u$ then 
\[ \mbox{cor}\left(z_x,z_{x'}\right)= \frac{2^{1-\nu}}{\Gamma(\nu)}\left(\sqrt{2\nu}u/\eta\right)^\nu K_\nu\left(\sqrt{2\nu} u/\eta\right), \]
where $K_\nu(.)$ is a modified Bessel function of the third kind. Let $\theta_1=(\eta, \nu)$. 

The second model, M2, is a marginally transformed $\chi^2_1$ field.   We simulate  from the Gaussian random field model M1, with parameters $\theta_2$.  Values are then squared to give $\chi^2_1$ margins. Let $y_x$ be the value at location $x$.  We next marginally transform to N(0,1) using  $z_x=\Phi^{-1}\{F_1(y_x)\}$, where $\Phi(.)$ and $F_1(.)$ are the N(0,1) and $\chi^2_1$ cumulative distribution functions, respectively.

The third model, M3, is a marginally transformed $F_{3,3}$ field.  We begin by constructing two independent $\chi^2_3$ fields, each from the sum of three squared independent
Gaussian random fields  M1 with  Mat\'{e}rn  parameters $\theta_3$.  The ratio of these $\chi^2_3$ fields has $F_{3,3}$ distribution and we again marginally transform back to N(0,1).

All simulated fields are mean-corrected to help with the comparisons.  The purpose of the marginal transformations is to generate fields that are not Gaussian random fields but which have Gaussian N(0,1) marginal distributions.  This is so that we can assess the added value of topological methods over and above simple comparisons of marginal distributions.  Similarly, in our simulations, for each choice of correlation $\theta_1$ for 
model M1 we choose  $\theta_2$ and  $\theta_3$  so as to match as closely as we can the correlation of the final M2 and M3 fields to Mat\'{e}rn with parameters $\theta_1$. Details of this are presented in Supplementary Material.  Again, the point is to show that differences between fields are not explained by differences in  correlation functions.  We have no particular interest in M2 and M3 other than as easily generated non-Gaussian random fields to be used as comparators.

The previously-presented Fig. \ref{fig:intro1} gives in panel (a) an example simulation from M1 on a $60 \times 60$ lattice with $\theta=(5,1).$  Panels (b) and (c) provide two examples of M3, and panel (d) an example of M2, with matched correlations.  The topological Nelson-Aalen plots for these fields are in Fig. \ref{fig:intro2}.  The solid line is the limiting value \eqref{eqn:limit}, which we obtained using the multivariate Gaussian distribution function routine in the R package {\tt mvtnorm}.

Table \ref{tab:cov1}   shows coverage of nominal 95\% confidence intervals and  bands  based on Nelson-Aalen plots for connected components from Gaussian random fields M1 simulated on a 60$\times$60 lattice.  Pointwise coverage is evaluated  at selected percentiles of the at-risk distributions, ie  levels $t$ satisfying $E\{Y_x(t)\}=p$ for five different values of $p$. The simultaneous confidence bands cover the central 90\% of the at-risk distributions, discretised to a grid of 200 equally spaced points. For each parameter combination we give results for three methods.  The first, for reference, is a naive approach that assumes  standard methods based on the variance estimator (\ref{eqn:var}) can be applied and using \cite{borgan87} for simultaneous confidence bands.   The second method assumes we have
a single replication and uses a parametric bootstrap with  $B=200$ simulated Gaussian random fields, generated  with  Mat\'{e}rn correlation parameters taken as the maximum likelihood estimates from the sole replicate.
The third method assumes we have replicates.   The standard and parametric bootstrap results are based on 1000 single simulations, the results with replications are from 1000 batches of $N=40$ replicates, which matches an application in the next subsection. 
There is severe undercoverage at times 
for the standard confidence intervals, but the other methods both produce good results.

Simulation results for M2 and M3 are provided in Supplementary Material,  along with further simulations to assess size and power when using our method to test for  a Gaussian random field.

\begin{table}
\begin{center}
\caption{Coverage of nominal 95\% pointwise confidence intervals and 95\% simultaneous confidence band (SCB) around Nelson-Aalen component plots for Gaussian random fields. Results from  1000 simulations  on $60 \times 60$ lattices,  with Mat\'{e}rn correlation function with parameters $\eta$ and $\nu$ }
\label{tab:cov1}  
\begin{tabular}{clrrrrrr}
 & &  \multicolumn{5}{c}{Percentile}\\
$(\eta, \nu)$  & Method   & 0.9 & 0.7 & 0.5 & 0.3 & 0.1 & SCB \\ 
\\
(5,1) & Standard        & 84.0& 92.2& 93.5& 88.9& 82.2  & 73.9\\
&Par. bootstrap& 95.3 &94.4 & 95.8& 95.3& 96.2& 94.1\\
&Replications    &  94.2 &94.7 &95.7 &93.9 &94.7 &94.4\\

\\ 
(10,1) & Standard   & 67.8 & 82.7 & 83.3 & 73.7 & 65.5 & 39.1\\
&Par. bootstrap&  95.3 &95.1 & 95.6 & 95.7 & 94.7 & 93.6 \\
&Replications    & 95.1 &94.0& 94.9 & 95.9 & 94.0 & 96.4\\
 \\

\\
(5,2) &  Standard  &    93.4 & 94.6 & 95.0 & 93.2 & 90.8 & 82.9\\     
&Par. bootstrap  & 95.4 &95.2& 95.7& 95.5 &95.8& 95.9\\
&Replications     & 94.0 & 94.6 &95.0& 94.3& 93.7 & 93.4\\
\end{tabular}
\end{center}
\end{table}

In Table \ref{tab:class} we provide a brief comparison between using the proposed Nelson-Aalen estimator and  existing topological data analysis techniques when discriminating between random fields, all on 60$\times$60 lattices.  The statistical properties of topological summaries such as persistence diagrams and landscape functions are not understood, which means we cannot use these methods to assess the adequacy of an assumed model such as a Gaussian random field, or for single-sample problems.  However, if we have two or more groups of replicated data  we can compare classification performance.  

We took  training data to be batches of 40 random fields simulated from models M1, M2 and M3, for each of the three parameter choices.   We then used a support vector machine to develop pairwise probabilistic classifiers, using the default parameter choices for the {\tt svm} routine in the R package {\tt e1071}.   We built classifiers based on the proposed Nelson-Aalen  estimator  $\hat{A}(t)$ defined at (\ref{eqn:nelaal})  and on the first 10 landscape functions \citep{bubenik2015} for each of components and holes, choosing the best of these to compare with $\hat{A}(t)$.   Landscape functions are popular summaries of persistence diagrams and are available through the R package {\tt TDA}.  We replicated the training data 20 times and assessed classifier performance on separate groups of 200 simulated random fields for each model.  

We assume test data would be allocated to the group with the highest class probability. In Table \ref{tab:class} we provide the accuracy of the classifier, which is the proportion of test data that are allocated to the correct group, and as a calibration measure we give the mean class probabilities for the correct groups.
 For example, a classifier that estimates the probability of the correct group to be 0.51  gives accurate allocation but is poorly calibrated, whereas one that estimates the correct class probability to be 0.99 would be both accurate and well calibrated.  The table  shows that classification based on the Nelson-Aalen estimator consistently outperforms classification based on the landscape functions in both accuracy and calibration.  Additional simulation results provided in Supplementary Material support this conclusion.

\begin{table}
\begin{center}
\caption{Classification of random fields.    Accuracy (Acc.) measures the proportion of random fields that were correctly allocated and calibration (Cal.) gives the mean estimated class probabilities for the correct classes. Values are multiplied by 100 for presentation}
\label{tab:class}  
\begin{tabular}{clrrrrrr}
Mat\'{e}rn &Method &  \multicolumn{2}{c}{M1 v M2}&  \multicolumn{2}{c}{M1 v M3}&  \multicolumn{2}{c}{M2 v M3}\\
 $(\eta, \nu)$  &   &  Acc.  & Cal.  &  Acc.   & Cal.  &  Acc.  & Cal.  \\
\\
(5,1) & Nelson-Aalen & 100.0 &97.7 & 99.4 & 96.6 & 100.0 & 97.5\\
 & Landscape      &96.5& 92.5& 91.7& 85.5&  95.4& 91.6\\
\\
(10,1) & Nelson-Aalen &99.7 &97.2 &95.2 &91.0 &100.0 &97.3\\
 & Landscape &99.7& 97.0& 80.7& 72.0&  95.8& 90.1\\
\\
(5,2) & Nelson-Aalen &100.0 &97.6 &99.7 &97.1 &100.0& 97.4\\
 & Landscape &  92.3& 87.7& 79.4& 71.0&  87.5 &80.1\\
\end{tabular}
\end{center}
\end{table}

\subsection{Applications: climate modelling, interstellar medium}

The Community Earth System Model Large Ensemble project \citep{kay2015community} is a publicly available ensemble of 40 high resolution climate model simulations. Although the models used to generate the data are deterministic, spatial and spatio-temporal statistical modelling is often applied as a means to both explore and concisely summarise the output \citep{cast16assess, edwards19}. Such models often assume that the modelled residuals form a  Gaussian random field. We can use the topological Nelson-Aalen approach to  assess this assumption, at least in part, and as a  quick and easy method to compare ensemble members and to examine the suitability of assumed correlation structures.

To illustrate our method we considered annual wind intensities at $10$m for the year 2020.  Each ensemble member has data on a 192$\times$288 grid over the surface of a sphere and at each data point we standardised by the mean and standard deviation across members.  To keep the correlation matrices manageable for this illustration we used
every third lattice point in both the latitudinal and longitudinal directions.  Following  
\cite{castruccio2014beyond} and \cite{castruccio2016compressing} we also trimmed the regions near the poles, leaving data on a $51 \times 96$ lattice for each ensemble member. 

\begin{figure}
\centerline{\includegraphics{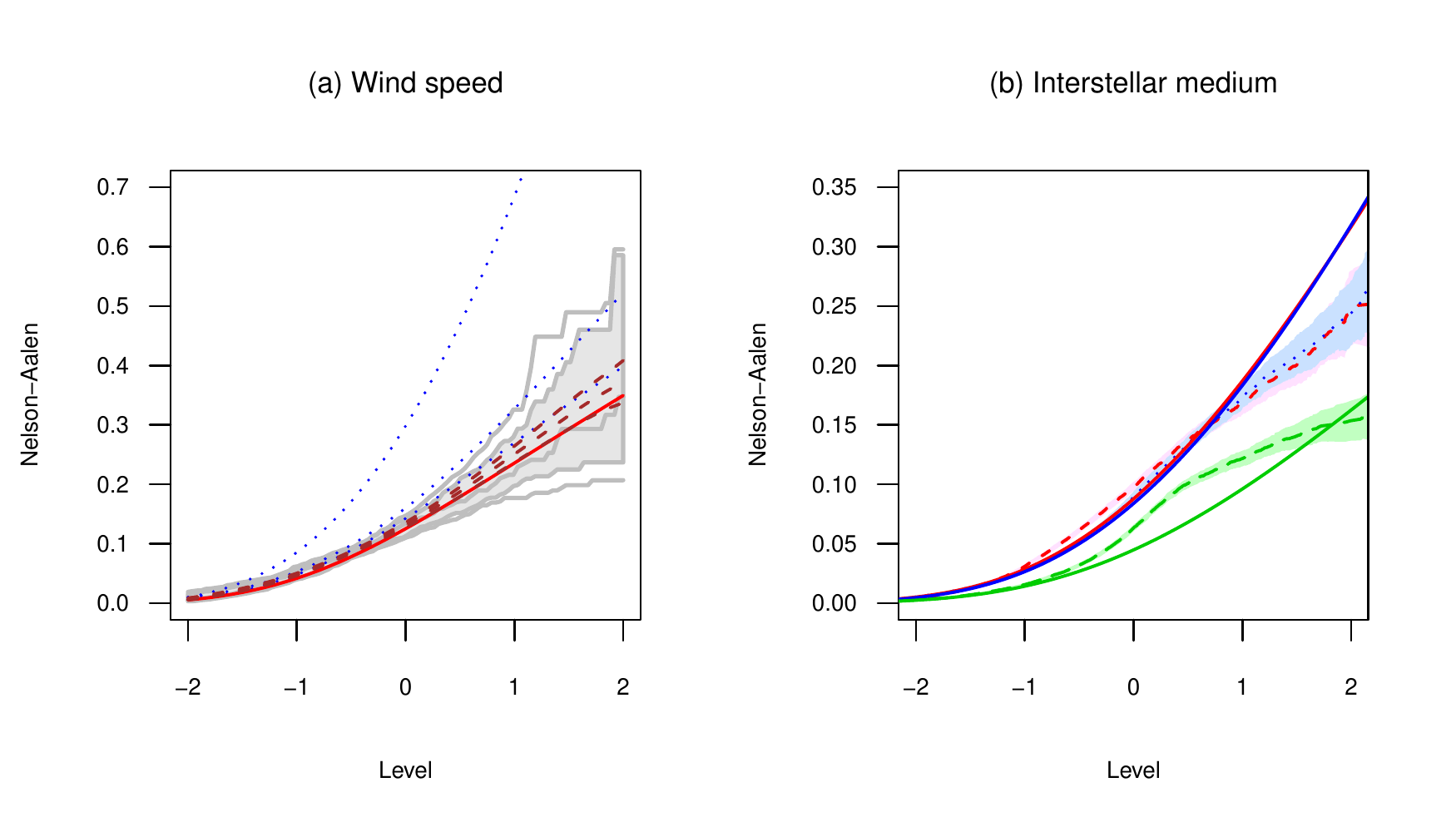}}
\caption{Nelson-Aalen plots for connected components. In (a) the grey region and lines represent functional boxplots of 40 separate realisations of climate data residuals.  The long dashed lines represent the mean and associated 95\% simultaneous confidence band, and the other lines are limiting values for Gaussian random fields with different correlation matrices. The solid line uses the empirical correlation matrix, the others are explained in the text.  Panel (b) considers the neutral hydrogen data: the short dashed lines are Region 1, the dotted lines Region 2 and the long dashed lines Region 3.  The shaded regions indicate $\pm$ two standard deviations.  The solid lines 
show the limiting values for Gaussian random fields with matched correlations.  The upper lines, which are almost indistinguishable, are for Regions 1 and 2, the lower is for Region 3. }
\label{fig:naexamples}
\end{figure}

Topological Nelson-Aalen plots for connected components are shown as functional boxplots \citep{sun11} in the left panel of  Fig. \ref{fig:naexamples} together with a simultaneous confidence band for the mean.  There is reassurance that the standard assumption of a Gaussian random field for these climate residuals is reasonable: the band almost uniformly includes the expected value for a Gaussian random field
 with correlation function taken to be the empirical correlation found in the data.  To illustrate the potential of our method as a diagnostic for fitted correlation models, the figure also includes expected values under three  modelled correlation structures. These are the three dotted lines in the plot. In order from the top,  the first line corresponds to a stationary 
Mat\'{e}rn model. It is clearly poor for these data, and  the same is true for the second dotted line, which corresponds to  a regional block model with separate Mat\'ern correlation functions fitted to land, ocean and coastal regions.  In contrast, as seen by the lowest dotted line, a semiparametric non-stationary model proposed by 
\cite{konzen2019} fits these data well.

Our second application is into the distribution in the Milky Way of neutral atomic hydrogen, \textup{H}\,{\sc i}, using data obtained  by the Galactic All-Sky Survey  using the Parkes 64 m radio telescope
\citep{mclure09, kalberla10}. Neutral atomic hydrogen makes up around 90\% of the atoms in the interstellar medium, which embeds stars and galaxies.  Turbulence and shocks lead to a highly heterogeneous random pattern to the local density of \textup{H}\,{\sc i}, and a common assumption in astrophysics is  that after the removal of trends and large-scale features, the distribution of the log-density of  \textup{H}\,{\sc i} is that of a Gaussian random field \citep{elmeg04, monin07}.  This assumption was investigated by \cite{hend2020}, who considered topological features of three distinct regions of the southern sky.  Based on counts of components and holes, but not event history, the authors concluded that two regions were consistent, after monotonic transformation, with Gaussian random fields but the third was not.

The right panel of Fig. \ref{fig:naexamples} shows Nelson-Aalen plots for connected components for these data.  Each region consists of density measurements over a 256$\times$256 lattice and following   \cite{hend2020} we considered the residuals after trend removal using a thin plate smoothing spline. The data were again marginally transformed to N(0,1) to preempt questions as to whether simple monotonic transformations could suffice to convert the data to Gaussian random fields.  Region 3 is clearly distinct from the others and not consistent with the fitted Gaussian random field. Regions 1 and 2  have very similar Nelson-Aalen plots and almost indistinguishable limiting values for Gaussian random fields with Mat\'ern correlations fitted to the data.  There is evidence however that at the higher levels these regions are not consistent with the fitted Gaussian random fields.

\section{COX PROPORTIONAL HAZARDS FOR EMBEDDED METRIC TREES}
\subsection{Topology and event history for tree data}

We now turn to trees $\Tcal=\{\Ncal,\Ecal\}$ made up of nodes $\Ncal$  and associated edges $\Ecal$ embedded  in $\mathbb{R}^2$ or $\mathbb{R}^3$, and with a single common ancestor.  A simple example was shown in  panel (a) of Fig. \ref{fig:trees}.

A number of approaches for the comparison of tree structures have been inspired by persistent homology. \cite{kanari2018} considered persistence of one-dimensional features on contours equidistant from the root, whether path length or radial distance, as previously illustrated in Panel (a) of Fig. \ref{fig:trees}. \cite{LI2017}  independently proposed a similar approach in the comparison of neuronal tree structures. \cite{BENDICH2016} 
applied persistent homology methods to brain artery trees in  $\mathbb{R}^3$,  \cite{MAOLI2017} used a geodesic distance function to derive persistence diagrams for the branching architectures of plants, and 
 \cite{BRODZKI2018}  used persistent homology to classify bronchial trees.

None of the methods of which we are aware can deal with the possibility of censored or incomplete trees.  For example,  an image may capture only part of a tree, in which case edges that cross the image boundary would  not be fully observed.  Edges that reach the boundary are not true first order nodes or leaves, and  should
not be considered as birth or death locations when applying any of the persistent homological techniques.  Equally they do contain information and these edges should not be ignored or discarded.

This type of problem is of course fundamental to event history and survival analysis methodology, where methods for censored or truncated data are very well established.    Hence our proposal is straightforward:  construct a filtration as above  but consider as an event the occurrence of either a leaf or a branching point.  More specifically we propose the following.

\begin{enumerate}
\item Filter outwards by radial distance $r$.
\item Consider a leaf event to occur at $r$ if a known  leaf has that radial distance from the root.  Consider a branch event similarly. If desired, we might distinguish bifurcations from trifurcations and so on.  Censored points cannot be events.
\item Construct a risk set consisting of all edges in $\Ecal$ that intersect a circle of radius $r$ centered at the root, as in panel (a) of Fig. \ref{fig:trees}.
\item Use standard event history methods to contrast events with risk sets.
\end{enumerate}

Let $\Ccal_r$ be  a disc of radius $r$ centred on the root and define $\Tcal_r=\Tcal \cap \Ccal_r$. Covariates can be included provided they are exogenous or dynamically defined on 
$\Ccal_r$.  Thus spatial location, branch order, edge weights, proximity to other branches in $\Ccal_r$ and so on can be handled easily using a Cox or other event time regression model. Further, we can construct risk sets from multiple trees and use fixed or random effects (frailties) to include tree-to-tree comparisons.

For this method to work we need three common event history assumptions.  First, that  events occur independently given the risk sets and covariates. Second, that censoring is independent.  Third, that whichever hazard or intensity model is assumed is appropriate. 
Under these assumptions standard results apply, including asymptotics for Cox or other model estimates.

\subsection{Application: vascular patterns in the eye}

\begin{figure}
\centerline{\includegraphics{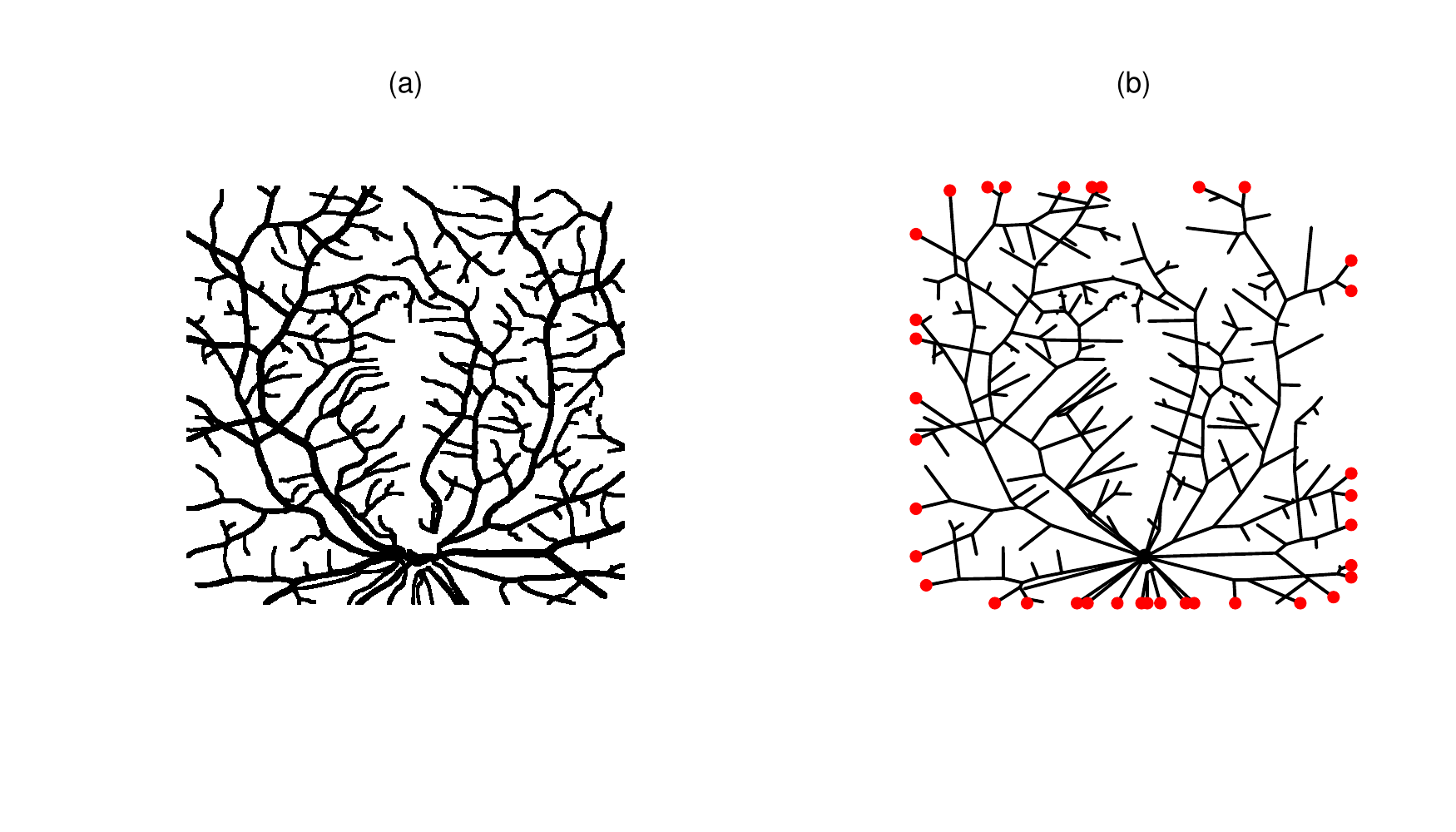}}
\caption{Left panel:  vascular patterns in a healthy human eye, with data taken from 
\cite{budai13}. Right panel: tree representation.  Censored edges are marked by points. }
\label{fig:eyefig}
\end{figure}

The left panel of Fig. \ref{fig:eyefig} shows the vascular pattern of a human eye, taken 
from a high dimensional image provided by \cite{budai13}.  The right panel of the figure shows a tree representation, with vertices joined by linear edges, a simplification of the optic nerve area to provide a root, and with the omission of  vessels that appear to originate outwith the plot region.  Vessels that reach the boundary of the plot region are marked as censored points.  This is a two-dimensional projection of a three-dimensional tree and branches can cross.  There are no cycles.

The data in Fig. \ref{fig:eyefig}  are from  one of the 15 fundus images of eyes of healthy patients provided by  \cite{budai13}, who also give fundus images of 15 patients with diabetic retinopathy.  We use our methods to compare diabetic retinopathy and healthy patients, and also to investigate the effect of a number of covariates.  These are the width of the vessel, the Euclidean distance between nodes, the ratio of path length to Euclidean distance between nodes as a measure of curvature, the number of other nodes within radius 200 units of the current node (about 10\% of the image size), the generation order, and the azimuth to the root measured from the vertical. We assumed log-linear effects for all covariates except azimuth, for which we assumed a  smooth effect via a $P$-spline with four degrees of freedom.    We also included patient-level terms as either fixed effects or frailties. In total over the 30 patients there were 12590 edges.  At this large sample size the fixed effect and frailty model fits were almost identical.

Table \ref{tab:cox} summarises some of the results.  In the upper part we consider an event to be the termination of a vessel at a leaf, and in the lower part we consider branching points.  Given the large sample, statistical significance is almost irrelevant and so to compare the magnitudes of effects we take the hazard ratio (hr) between 
the 20\% and 80\% points of the  covariate distributions.  

As might be expected, large vessels are highly unlikely to terminate in a leaf but are more likely to branch.  The hazard functions for both types of event fall as Euclidean distance from an earlier node increases.  There is no real effect of either curvature of an edge or proximity to other nodes, but there is some evidence that branching becomes less likely as the generation of an edge increases.  There is a modest azimuth effect on branching, with a plot of the smooth fit (not shown here) indicating a decrease in branching hazard as aziumth moves away from the vertical, fairly symmetrically in both directions.

\begin{figure}
\centerline{\includegraphics[height=3in,width=3in]{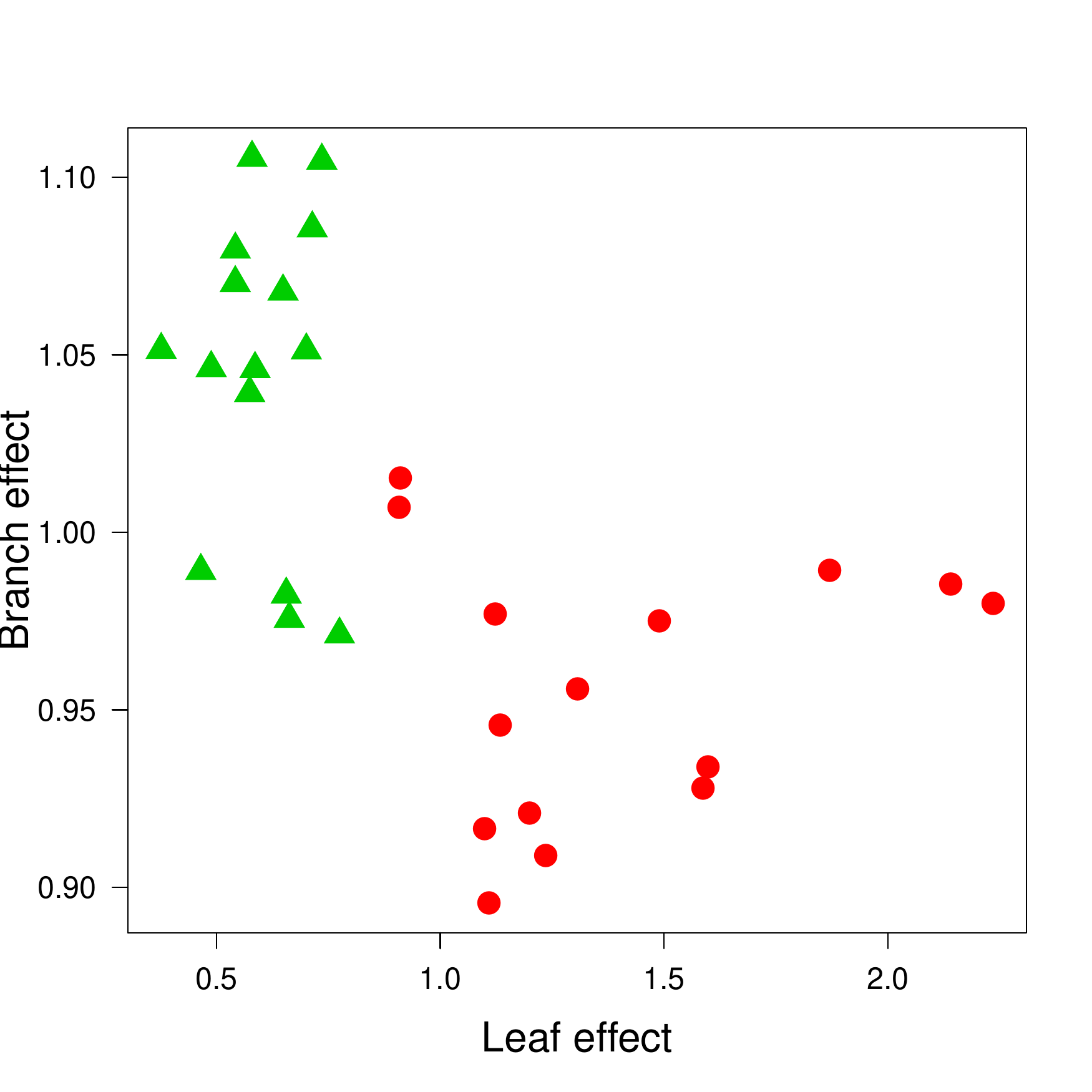}}
\caption{Subject-level frailties for 15 diabetic retinopathy patients (triangles) and 15 healthy patients (circles).}
\label{fig:frailfig}
\end{figure}

Of most interest is the difference between patients.  The estimated frailty terms for both leaf and branch events are shown in Fig. \ref{fig:frailfig}.  There are relatively small differences between patients in vessel branching frequency, but large differences in the hazard for terminating in a leaf, given covariate effects. Interestingly, we have complete separation between healthy and diabetic  retinopathy patients in these frailties.

\begin{table}[h]
\begin{center}
\caption{Cox proportional hazards results for eye trees}
\label{tab:cox}
\begin{tabular}{lrrrrr} 
 & \multicolumn{5}{c}{Event = leaf}\\
  & Coefficient & SE & $\chi^2$ & df & hr\\ 
Width & -0.3849 & 0.0079 & 2369.51 & 1 & 0.02\\ 
Euclid & -0.0051 & 0.0002 & 797.16 & 1 & 0.41\\ 
Path length ratio & 1.0142 & 0.1260 & 64.76 & 1 & 1.05\\ 
Nodes within 200 & 0.0037 & 0.0037 & 1.02 & 1 & 1.02\\ 
Order & -0.0147 & 0.0049 & 8.90 & 1 & 0.86\\ 
Azimuth &  &  & 34.66 & 4 & 1.18\\ 
Patient effect &  &  & 823.98 & 29 & 2.29\\ 
\\
 & \multicolumn{5}{c}{Event = branch}\\
 & Coefficient & SE & $\chi^2$ & df & hr\\  
Width & 0.0591 & 0.0016 & 1420.28 & 1 & 1.88\\ 
Euclid & -0.0058 & 0.0001 & 1874.37 & 1 & 0.37\\ 
Path length ratio & -0.9058 & 0.2406 & 14.18 & 1 & 0.95\\ 
Nodes within 200 & 0.0375 & 0.0033 & 128.54 & 1 & 1.21\\ 
Order & -0.0622 & 0.0051 & 149.35 & 1 & 0.54\\ 
Azimuth &  &  & 119.30 & 4 & 1.41\\ 
Patient effect &  &  & 103.60 & 29 & 1.37\\ 
\end{tabular}
\end{center}
\end{table}

\section{DISCUSSION}

There are a number of areas for further work.  In Section 3 we suggested empirical variances for the Nelson-Aalen estimator when there are repeated samples of random fields, as in the climate example, and parametric bootstrap estimators when there was interest in a particular model, such as a Gaussian random field. In Supplementary Material,  an adjusted parametric bootstrap assuming local misspecification as in \cite{Copas2005} is suggested and
provides a conservative approach for the case of a single sample and no model. However, it  would be useful formally to develop a non-parametric variance estimator. 

In both Section 3 and Section 4 we concentrated on the birth times of connected components. We might also consider death times, though the identification of risk sets is less straightforward in this case, which in turn means that limiting properties will be harder to derive.  For complexes in two or higher dimensions, including the random fields of Section 3, we could consider higher order features.  In two-dimensions it will be particularly simple to use our methods to investigate death times of holes, which as a result of Alexander duality \citep{edel12} can be analysed in the same manner as births of components simply by reversing the direction of the filtration. 

\section{ACKNOWLEDGEMENT}

We thank Evandro Konzen and Jian Shi for providing software to fit their non-stationary correlation model to the climate data. 



\clearpage

\bibliographystyle{biometrika}
\bibliography{toprefs}

\end{document}